\input amstex
\input epsf
\documentstyle{amsppt}

\def\esim{\underset{e}\to\sim}
\def\wesim{\underset{e}\to\asymp}
\def\ewsim{\wesim}

\def\RR{\Bbb{R}}

\define\eps{\varepsilon}
\define\sign{\operatorname{sign}}
%
\define\lempnOne   {10}
\define\lempnTwo   {11}
\define\lempnThree {12}
%
%

\define\figOvalsTree    {1}
\define\figEnumerLeaves {2}
\define\figRemoveEdge   {3}


\def\refArn     {1}
\def\refDG      {2}
\def\refHP      {3}
\def\refKnuth   {4}
\def\refKO  {5}
\def\refOtter   {6}
\def\refPolya   {7}

\def\refRS  {8}
\def\refS       {9}


\def\eqOtter  {1} 
\def\eqConst  {2} 
\def\eqKO     {3} 
\def\eqAr     {4}  
\def\eqLbound {5}  
\def\eqhalf     {6} 
\def\eqT      {7} 
\def\eqFE     {8} 
\def\eqRR     {9} 
\def\eqFy     {10}
\def\eqa      {11}
\def\eqAsOne    {12}
\def\eqrprime   {13} 
\def\eqrsecond  {14}

\rightheadtext{
Trees half of whose vertices are leaves and 
curves enumeration}
\topmatter
\title    
	The number of trees half of whose vertices are leaves
	and asymptotic enumeration of plane real algebraic curves
\endtitle
\author
        V.M.~Kharlamov, S.Yu.~Orevkov
\endauthor

\address
	First author: Universit\'e Louis Pasteur et IRMA (CNRS),
7 rue Ren\'e Descartes 67084 Strasbourg Cedex, France
\endaddress

\email
	kharlam\@math.u-strasbg.fr
\endemail

\address
	Second author: Laboratoire Emile Picard, UFR MIG, Univ. Paul Sabatier,
	118 route de Narbonne, 31062 Toulouse, France
\endaddress

\email
	orevkov\@picard.ups-tlse.fr
\endemail

\abstract 
The number of topologically different
plane real algebraic curves of a given degree $d$ has the form
$\exp(C d^2 + o(d^2))$. We determine the best available upper 
bound for the constant $C$.
This bound follows from Arnold inequalities on the number of empty ovals. To
evaluate its rate we show its equivalence with the rate of growth
of the number of trees half of whose vertices are leaves and
evaluate the latter rate.
\endabstract

\thanks
First author is a member of Research Training Networks EDGE and RAAG,
supported by the European Human Potential Program.
\endthanks

\endtopmatter

\document

\subhead Introduction
\endsubhead
Recall that a {\it rooted tree} is a tree with a distinguished
vertex. The distinguished vertex is called the {\it root}. The
{\it multiplicity} or the {\it valence} of a vertex is the number
of edges which are incident to it. A vertex of multiplicity one is
called a {\it leaf}. By convention, we assume that the root is a leaf
if the tree has no other vertices. Otherwise, the root is not
considered as a leaf even if its multiplicity is one.

In this paper we work exclusively with {\it unlabelled finite
trees}.

Rooted unlabelled trees are used to encode the topology of
nonsingular
curves in the real projective plane (by a nonsingular curve we
mean a closed one-dimensional, not necessarily connected,
sub-manifold).
We associate the vertices
with the connected components of the complement of the curve. The
root will correspond to the component with non-oriented closure
and the tree will represent the adjacency relations between the
components (see Figure \figOvalsTree). The fact that this graph is
a tree follows from the Jordan curve theorem. It is finite since
our curves are compact.

It is worth noticing that two curves have the same encoding if and
only if there is an ambient isotopy transforming one into
another.

\midinsert
\epsfxsize 80mm
\centerline{\epsfbox{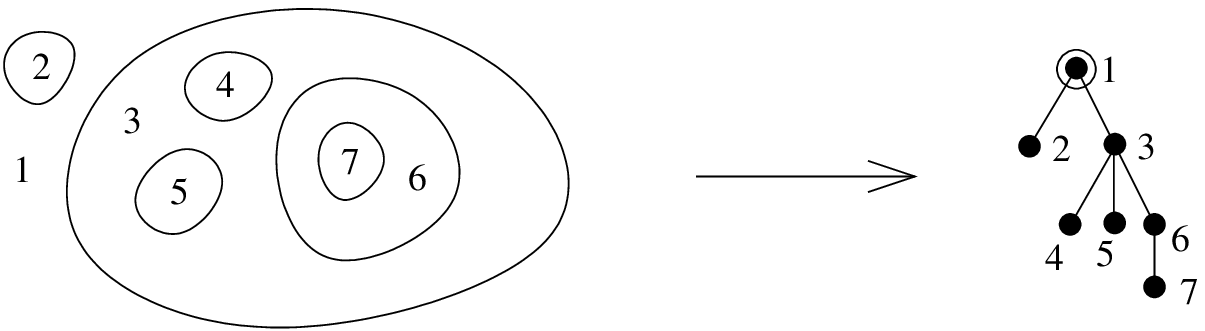}}
\botcaption{ Figure \figOvalsTree. Plane curve and the corresponding rooted tree }
\endcaption
\endinsert

Even if the curves are algebraic, there is no
any restriction on the
encoding tree as long as no condition on the curve is imposed. The
situation is changing as soon as we fix the degree $d$ of the
curve. Then, already the number of
connected components, and thus the number of the vertices in the
encoding tree, is not arbitrary.
According to our convention, it is the two-sided components of the
curve which represent the edges, so that the number of vertices is
the same as the number of components of the curve
if $d$ is odd, and it is greater by $1$
if $d$ is even. As is
known, the number of components of the curve is $\le
\frac{(d-1)(d-2)}2 +1$. Introduce, thus, the following
notation
$$
    N_d = \cases (d-1)(d-2)/2+1, &\text{if $d$ is odd,}\\
        (d-1)(d-2)/2+2, &\text{if $d$ is even.}\endcases
$$

Starting from $d=4$, not any tree with $\le
N_d$ vertices can be realized by a curve of degree $d$. Let $I_d$
be the number of the trees which can be realized by curves of
degree $d$. No direct formula or functional equation for
these numbers is known;
moreover,
their exact values
are available only for $d\le 7$. Very few is known even on the
rate of growth of $I_d$.

As is shown in \cite{\refKO},
$$
I_d \ewsim \exp(d^2),
$$
where $a_m\ewsim b_m$ means that
$\log a_m=O(\log b_m)$ and $\log b_m=O(\log a_m)$.
On the other hand,
due to Otter \cite{\refOtter} (see also \cite{\refHP; Section
9.5}), one has the following exponential equivalence for the
number $T_n$ of rooted unlabelled trees with $n$ vertices
$$
        T_n\esim C^n,
         \qquad C=2.95576\dots,                           \eqno(\eqOtter)
$$
where the latter means that $\log T_n\sim n\log C$. This implies
that
$$
T_1+\dots+T_n\esim C^n,
$$
hence,
$$
    I_d \le C^{{d^2\over2} + o(d^2)}.
                                                 \eqno(\eqConst)
$$

The aim of the present note is to correct one erroneous remark
from \cite{\refKO} and to show that the so-called Arnold
inequalities \cite{\refArn} allow to reduce the constant $C$ in
the estimate (\eqConst). Namely, we prove that according to these
inequalities
$$
I_d \le C_1^{{d^2\over2}+o(d^2)},\quad
C_1=2.9193800\dots \eqno(\eqKO)
$$
More precisely,
{\it $(\log C_1) d^2/2$ is asymptotically equivalent to $\log A_d$
where $A_d$
is the number of unlabelled trees with $n\le N_d$ vertices not
excluded by the Arnold inequalities}. It shows that the Arnold
inequalities exclude more arrangements of $\le N_d$ closed simple
circuits than any other known property of plane algebraic curves,
including (see \cite{\refKO})
 the consequences of the Bezout theorem.

Let us recall that the principal
Arnold inequalities concern the curves of even degree
$d=2k$ exclusively. They state that
$$
{\text{even}}^*\le \frac{(k-1)(k-2)}2, \quad
{\text{odd}}^*\le \frac{(k-1)(k-2)}2,
\eqno(\eqAr)
$$
where
$\text{even}^*$ is the number of not end vertices of odd
distance from the root, and
$\text{odd}^*$ is the number of not end vertices of even (non
zero) distance from the root. These inequalities imply the
following  lower bounds on the number $l$ of leaves (end vertices)
whatever is the parity of the degree $d$:
$$
l\ge n-1-\left[\frac{d-1}2\right]\left(\left[\frac{d-1}2\right]-1\right)
\eqno(\eqLbound)
$$
where $n$ is the total number of vertices. If $d$ is even it is a
straightforward consequence of (\eqAr) and if $d$ is odd it
follows from (\eqLbound) for $d+1$. In particular, for the maximal
value $n=N_d$ of $n$, the right hand side is approximately the
half of $n$:
$$
	N_d-1-\left[\frac{d-1}2\right]
	\left(\left[\frac{d-1}2\right]-1\right)
	\sim\frac12 N_d.
							\eqno(\eqhalf)
$$
According
to results of this note, it is
the arrangements  with $n=N_d$ and $l\sim \frac12 N_d$ which
determine the asymptotical impact of Arnold bounds: {\it $A_d$ has
the same $\ewsim$-rate of growth as the number of the trees with
$N_d$ vertices half of which are leaves.}
In particular, the upper bound for $I_{2k}$ deduced from the sole inequality (\eqLbound)
has the same $\ewsim$-rate of growth as the upper
bound which can be deduced from (\eqAr).

In fact,
what is important in the coefficient $1/2$ in
(\eqhalf)
is that $1/2>0.438156...$
If the Arnold inequalities
were not known but someone
proved
only that
$l > 0.43 N_d$, this fact would not reduce the constant $C$ in (\eqConst)
because the most of trees have about $43.8\%$ leaves (see Appendix for details
and references).

The note is organized as follows. The asymptotic growth of the
number of the trees half of whose vertices are leaves is
established in Section 1 in Theorem
7. The asymptotic impact of the Arnold inequalities is deduced
from this theorem in Section 2: 
Theorem 9 takes into account only the
bound (\eqLbound) and Theorem 13 shows that (\eqAr) does not
improve the rate. In Appendix we compare the result with the
limiting distribution and show that the central limit theorem is
not sufficient for our purpose: the range of values we treat is
outside the range of a suitably good convergence.

\head  1. On trees half of whose vertices are leaves.
\endhead

\subhead 1.1. Functional equation
\endsubhead
Let us denote the number of rooted unlabelled trees with $n$
vertices and $k$ leaves by $a_{n,k}$
and consider the associated bi-variant generating function (a formal power
series)
$$
    T(x,z) = \sum_{n,k} a_{n,k} x^n z^k = \sum_{n=1}^\infty a_n(z)
    x^n.
                                \eqno(\eqT)
$$
We get (see Figure \figEnumerLeaves)
$$
    T(x,z) = zx + zx^2 + (z + z^2)x^3 + (z + 2z^2 + z^3)x^4 +
         (z + 4z^2 + 3z^3 + z^4)x^5 +
$$
$$
    (z+6z^2+8z^3+4z^4+z^5)x^6 + (z+9z^2+18z^3+14z^4+5z^5+z^6)x^7+\dots
$$

\midinsert
\epsfxsize 120mm
\centerline{\epsfbox{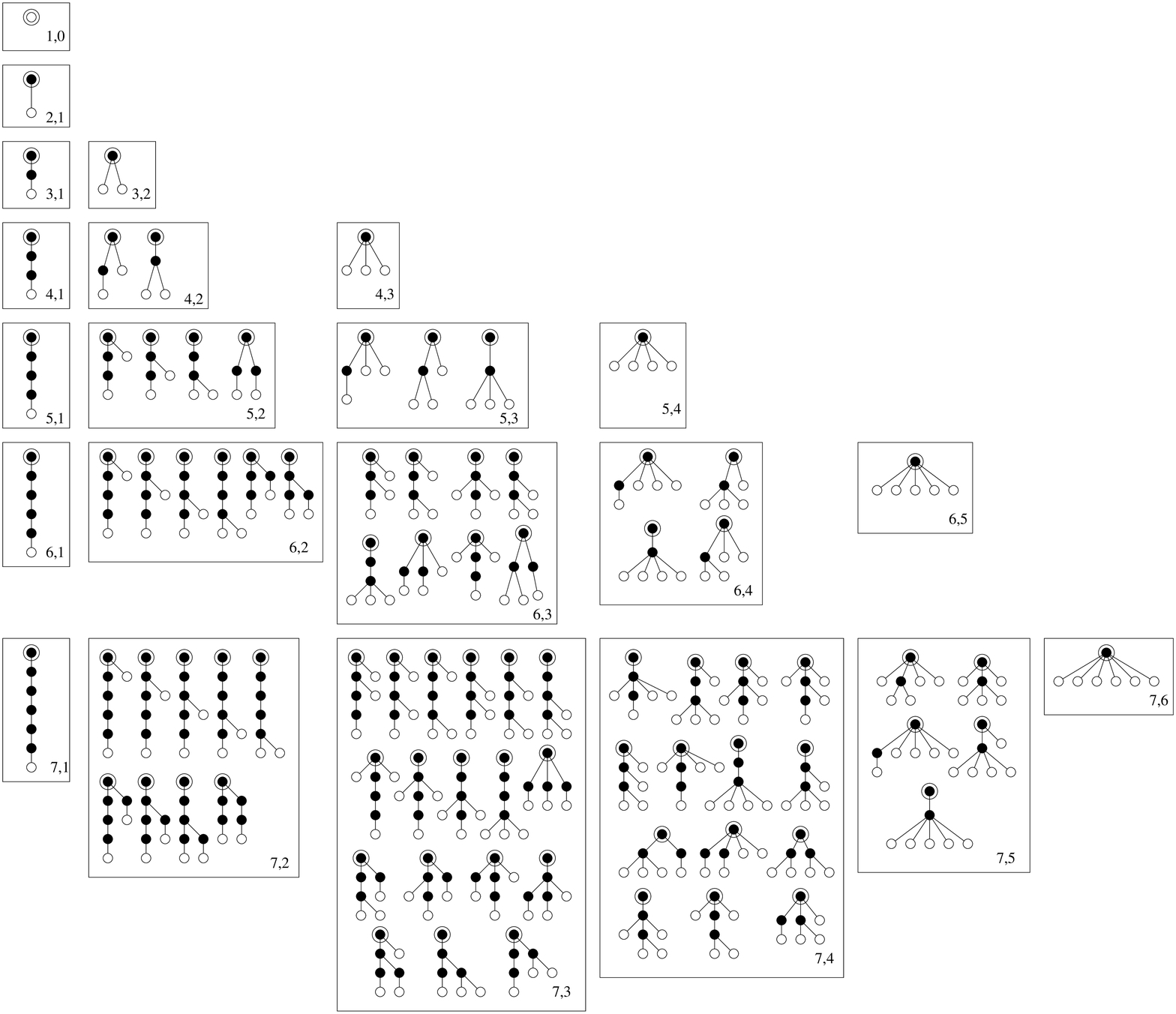}}
\botcaption{ Figure \figEnumerLeaves }
\endcaption
\endinsert

For technical reasons, we
introduce also
$$
    \tilde T(x,z) = T(x,z) - zx + x = \sum_{k=1}^\infty\tilde a_k(z)x^k,
    \qquad
    \tilde a_k(z) = \cases 1, &k=1\\ a_k(z), &k>1,\endcases
$$
which is the generating function under the
convention
that the vertex of the one-vertex tree is not considered
as a leaf.

Using P\'olya enumeration theorem as it is done in [\refRS] one
can prove that $T(x,z)$ satisfies the (formal) functional equation
$$
    {\tilde T(x,z)} = {T(x,z)} - zx + x =
        x\exp\left(\sum_{k=1}^\infty {T(x^k,z^k)\over k}\right).
                                \eqno(\eqFE)
$$
The specialization $T(x) = T(x,1)$ is the classical generating
function for the number of rooted unlabelled trees and
substituting of $z=1$ into (\eqFE) turns it into the classical
P\'olya equation, see [\refPolya].

It may be worth noticing that to prove (\eqFE) one can use as well
the following bi-variant analog of the Cayley
product formula for $T(x)$, cf. [\refKnuth],
$$
\tilde T(x,z)=\frac{x}{\prod(1-x^nz^m)^{a_{n,m}}}.
$$


\subhead 1.2. Recurrent relation
\endsubhead
Taking the logarithmic derivatives of the both sides of (\eqFE), we get
$$
    {\tilde T_x(x,z)\over \tilde T(x,z)}
    = {\partial\over\partial x}\left(\log x
        + \sum_{k=1}^\infty {T(x^k,z^k)\over k}\right)
    = {1\over x} + \sum_{k=1}^\infty x^{k-1}T_x(x^k,z^k)
$$
Multiplying the both sides by $x\,\tilde T(x,z)$ and subtracting $\tilde T(x,z)$, this gives
$$
    x \tilde T_x(x,z) - \tilde T(x,z) = \tilde T(x,z)\sum_{k=1}^\infty x^k T_x(x^k,z^k).
$$
Hence,
$$
    \sum_{n=1}^\infty
        n \tilde a_{n+1}x^{n+1}
    = \sum_{p=1}^\infty\tilde a_p x^p \sum_{k=1}^\infty \sum_{j=1}^\infty ja_j(z^k)x^{jk}
    = \sum_{n=1}^\infty x^{n+1}\sum_{p+jk=n+1} ja_j(z^k)\tilde a_p(z).
$$
Thus, we obtain the
recurrence relation (cf. [\refOtter] and [\refRS])
$$
    na_{n+1}(z) =n\tilde a_{n+1}(z)
        = \sum_{j=1}^n j\sum_{k=1}^{[n/j]} a_j(z^k)\tilde
        a_{n+1-jk}(z).
                                        \eqno(\eqRR)
$$
Together with the initial conditions $a_1(z)=z$, $\tilde a_1(z)=1$, the relation
(\eqRR) gives a rather fast way to compute $a_n(z)$.


\subhead 1.3. Analytic properties of $T(x,z)$
\endsubhead
If before we treated the generating functions as
formal series, now we need to study their analytic behavior.

Let $\alpha$ be the radius of convergence
of the
power series $T(x)$.
Using
Polya's
approach, see [\refPolya], i.e., resolving the equation
$x\exp\big(1+\sum_{k=2}^\infty T(x^k)\big)=1$ (for instance, by
Newton's method), one can compute $\alpha$ with any given
precision. Indeed, any finite nimber of coefficients of the
involved series can be computed using (\eqRR) and the number of
terms to be summated, can be found from some rough estimate of
$\alpha$. Performing this computation, one gets
$$
    \alpha = 0.33832185689920769519611262571701705318...
$$
This constant is sometimes called Otter constant because the first seven digits
were computed in [\refOtter] (using the
above approach from [\refPolya]).

Let us denote by $D$ the domain of convergence of the series
(\eqT). Here, we follow the classical tradition and mean by the
{\it domain of convergence} the interior of the set where the series is
convergent. As is known, it coincides with the interior of the set
of points $(x,z)$ such that $\sup_{n,k}|a_{n,k}x^nz^k|<\infty$. An
important, also well known, consequence is that the logarithmic
image
$$
\log |D| = \{(\log|x|,\log|z|) \,:\, (x,z)\in D\}\subset \Bbb R^2
$$
of any
convergence domain is convex (in other words, the convergence
domains are {\it logarithmically convex}).

\proclaim{ Lemma 1} There exists a continuous function
$\zeta\mapsto r(\zeta), \RR_{>0}=\{\zeta>0\}\to\RR_{>0},$
such that $D = \{(x,z)\,:\, |x|<r(|z|)\}$.
Moreover, $\alpha/\zeta\le r(\zeta)\le \alpha$ for $\zeta\ge 1$ and
$r(\zeta)<\min\{1,\frac1{|\zeta|}\}$ for any $\zeta$.

The series $T(x,z)$ converges at each point $x=r(z), z>0,$
of $\partial D\cap\Bbb R^2_{>0}$.
\endproclaim

\demo{Proof} Due to logarithmic convexity of $D$, the existence
statement follows from $D\subset\{|xz|\le 1\}$; in its turn, this
inclusion follows from $a_{n,k}\ge 1$ for any $n>k$. The proof of
other non strict bounds uses, in addition, the cited above
convergency properties of $T(x)=T(x,1)$. The strict inequality
$r(z)<\frac1{|z|}$ is a consequence of the convergence of $T(x,z)$
at the boundary points. To prove this convergence it sufficient to
notice that
$$T(x,z) = xz-x+x e^{T(x,z)+\dots} > xz-x+x e^{T(x,z)}$$
for $0<x<r(z), z>0$; it implies the boundedness of $T$ on the interval $x\in[0,r(z)[$
and, by Abel theorem,
its convergence at $x=r(z)$.
\qed\enddemo

\proclaim{Lemma 2}
The transformations $(x,z)\mapsto(x^k,z^k), k\ge 2,$ map $D$ into itself.
For any point $(x,z), z\ne 0,$ in the closure of $D$ the series
$$
 h(x,z) = \sum_{k=2}^\infty {T(x^k,z^k)\over k}
$$
is absolutely convergent and defines a function holomorphic at
such a point.
\endproclaim
\demo{Proof} The invariance property follows from
the logarithmic convexity and the bounds on $r(z)$ given by Lemma 1.
In addition, due to this Lemma, for all $(x,z)$ in a small neighborhood
of any point in the closure of $D$ we have bounds $|x^k|\le a^k,
|z^k|\le b^k$ with $a<1, ab <1$ whatever is $k\ge 1$. These bounds provide
a bounded convergence of the series:
$$
\sum_{k\ge 2}\sum_{n,m}\frac{|a_{n,m}x^{nk}z^{mk}|}k\le
\sum_{n,m}\sum_{k\ge 2}\frac{a_{n,m}a^{nk}b^{mk}}k
\le\lambda \sum_{n,m}a_{n,m}a^{2n}b^{2m}=\lambda T(a^2,b^2).\qed
$$
\enddemo

In what follows we study the boundary values
$a(z)=T(r(z),z), z>0$ of $T$ and use an auxiliary function
$$
    F(x,y,z) = z - 1 + e^{y + h(x,z)} - \frac{y}x.
$$
By (\eqFE), we have
$F(x,T(x,z),z)=0$ at any point of the closure of $D$
with $x\ne 0, z\ne 0$.
In particular, the real curve $x=r(z), z>0,$ satisfies the
equation
$$
F(x, a(z), z)=0.
$$

\proclaim{ Lemma 3}
The function $r(\zeta)$ is analytic.
The function $F(x,y,z)$ is analytic near the real curve $x=r(z)$, $z>0$.
We have
$$
    F_y(r(z),a(z),z) = 0,           \eqno(\eqFy)
$$
$$
    a(z) = 1 + r(z)(z - 1).          \eqno(\eqa)
$$
\endproclaim

\demo{ Proof }
The analyticity of $F$ follows from Lemma 2, and then all the other
statements, except the relation (\eqa), follow from the implicit
function theorem.

Let show that $a(z) = 1 + r(z)(z - 1)$.
By the definition of $F$, we have
$F_y= e^{y+h(x,z)} - \frac1x$. Hence, for $x=r(z)$ and $y=a(z)$ we have
$$
    0=F_y = e^{y+h(x,z)} - \frac1x
    \qquad\text{and}\qquad
    0=F = z-1 + e^{y+h(x,z)} - \frac{y}x.
$$
Thus, $\frac{y}x = \frac1x + z - 1$ and $y = 1 + x(z-1)$.\qed
\enddemo

Due to Lemma 3, the
function $x=r(z)$ can be found by resolving the equation
$$
    x\exp(1 + (z-1)x + h(x,z)) = 1
$$
This allows one to compute $r(z)$ with any given precision.

Let us
define
$$
\split
    &a_n^+(z) = \sum_{k>n/2} a_{n,k} z^k,   \quad
    a_n^-(z) = \sum_{k\le n/2} a_{n,k} z^k, \quad
    T_\pm(x,z) = \sum_{n=1}^\infty a_n^\pm(z)x^n;
\\
    &\hat T(x,z) = \hat T(x z^{-1/2},z),    \quad
    \hat T_\pm(x,z) = \hat T_\pm(x z^{-1/2},z), \quad
    \hat r(\zeta) = r(\zeta)\sqrt{\zeta}.
\endsplit
$$
and denote by
$\hat D$ and $\hat D_\pm$
the domain of convergence of $\hat T$ and $\hat T_\pm$, respectively.
It is clear that $\hat D=\{(x,z)\,:\, |x|<\hat r(|z|)\}$.

\proclaim{Lemma 4}
The function $\hat r(\zeta)$ has a single critical point, this point
is a point of maximum.
\endproclaim
\demo{Proof}
The logarithmic map
$(x,z)\mapsto(\log|x|,\log|z|)$ transforms $xz^{-\frac12}$
in a linear function. Therefore, due to the convexity of $\log|D|$, the critical points
of $\hat r(\zeta)$ form a convex set. If it is not reduced to a single point, then,
since $r$ is real analytic, $r(\zeta)=c\zeta^{-\frac12}, c>0,$ which contradicts
to the bounds from Lemma 1.

It is a point of maximum, since the domains of convergence are Reinhardt domains,
i.e., $(x,z)\in D$ as soon as there exists $(x_0,z_0)\in D$ with $|x|<|x_0|, |z|<|z_0|$.
\qed\enddemo

Denote by $z_0\in \RR_+$
the point where the maximum of $\hat r(\zeta)$ is attained and put $x_0=\hat r(z_0)$.

\proclaim{Proposition 5}
$\hat D_\pm=\{(x,z)\,:\, |x|<\hat r_\pm(|z|)\}$, where
$$
    \hat r_-(\zeta) = \max_{\omega\le\zeta}\hat r(\omega)
    = \cases \hat r(\zeta), &\zeta\le z_0,\\ x_0, &\zeta\ge z_0, \endcases
    \,\,\,\text{and}\,\,\,\,
    \hat r_+(\zeta) = \max_{\omega\ge\zeta}\hat r(\omega)
    = \cases \hat r(\zeta), &\zeta\ge z_0,\\ x_0, &\zeta\le z_0. \endcases
$$
\endproclaim

\demo{ Proof }
For a point $p=(u_0,v_0)\in\RR^2$, let us denote
$\RR^2_{+-}(p) = \{(u,v)\,|\, u\le u_0, v\ge v_0\}$ and
$\RR^2_{++}(p) = \{(u,v)\,|\, u\le u_0, v\le v_0\}$.
The result follows from the following properties:

\smallskip
(a) $\log|\hat D|$ and $\log|\hat D_\pm|$ are convex;

(b) If $p\in
\log|\hat D_\pm|$ then $\RR^2_{+\pm}(p)\subset
\log|\hat D_\pm|$;

(c) 
$\log|\hat D| =\log|\hat D_+| \cap \log|\hat D_-|$.
\qed\enddemo

\subhead 1.4. Rate of growth
\endsubhead
\proclaim{Theorem 6}
$$
    \sum_{k>n/2} a_{n,k} \underset{e}\to\sim C_1^n,
    \qquad\text{where}\quad
    C_1 =
    \frac1{x_0} =
    2.919380017448416911265032583985...
$$
\endproclaim
\demo{Proof}
The coefficients $a_n^+(1)=\sum_{k>n/2} a_{n,k}$ of the power series
$\hat T_+(x,1)$ satisfy the following relation
$$
\log a_{n+m+2}^+(1) \ge \log a_{n}^+(1)+\log a_{m}^+(1)-\log 2
$$
(to prove this relation it is sufficient to plant two trees over a
new root and to add a leaf growing from the root). Hence, the
sequence
$n^{-1}\log a_n^+(1)$
has a limit and, by the Cauchy
rule,
$$
    \sum_{k>n/2} a_{n,k} = a_n^+(1) \underset{e}\to\sim \hat r_+(1)^{-n}.
                                \eqno(\eqAsOne)
$$

To compute $\hat r_+(1)$, we must find $z_0$. We compute it as the
root of the equation $\hat r'(z)=0$ (the root is unique by the 
convexity of $\log D$).
To find it by Newton's
method, we need $\hat r'(z)$ and $\hat r''(z)$. They
can be found as follows. Derivating the identity
$F(r(z),a(z),z)=0$ and using (\eqFy), we get
$$
    F_x(r(z),a(z),z)r' + F_z(r(z),a(z),z) = 0.  \eqno(\eqrprime)
$$
Derivating again, we see that at points $(r(z),a(z),z)$ one has
$$
    F_{xx}{r'}^2 + F_{xy}r'a' + 2F_{xz}x' +
    F_{yz}a' + F_{zz} + F_x r'' = 0.
                    \eqno(\eqrsecond)
$$
Note that $a'$ can be found from (\eqa).

The partial derivatives of $F$ at a point $(r(z),a(z),z)$ are
$$
    F_x=(h_x/r)+(a/r^2),\qquad
    F_y=0,\quad
    F_z=1+(h_z/r),
$$
$$
    F_{xx}=(h_{xx}+h_x^2)/r - 2(a/r^3),\qquad
    F_{xy}=(h_x/r)+(1/r^2),\qquad
    F_{xz}=(h_{xz}+h_x h_z)/r,
$$
$$
    F_{yz}=h_z/r,\qquad
    F_{zz}=(h_{zz}+h_z^2)/r
$$

Solving the equation $\hat r'(z)=0$ by Newton's method, we find
$$
    z_0 = 1.48491739577413809587489...
$$
and
$$
    x_0 = \hat r(z_0) = 0.3425384821514313844959919944869...
$$
Since $z_0>1$,
we have $\hat r_+(1) = \hat r(z_0) = x_0.$
Now, the desired asymptotic relation follows from
(\eqAsOne)
and
$$
   C_1 = 1/x_0 = 2.919380017448416911265032583985...\qed
$$
\enddemo
\proclaim{Theorem 7}
There is a continuous function $\lambda\mapsto C(\lambda)$,
$\RR_{\ge 0}\to\RR_{\ge 0}$, such that
$$
    \sum_{k>\lambda n} a_{n,k}\underset{e}\to\sim C(\lambda)^n
\quad\text{
for any $\lambda \ge 0$.}
$$
For each $\lambda>\frac12$ one has $C(\lambda)<C(\frac12)=C_1.$
\endproclaim
\demo{Proof} Let $z_{0,\lambda}$ be the critical point of
$r(\zeta)\zeta^\lambda$. By the same arguments as in the proof of
Proposition 5 and Theorem 6,
$$
    \sum_{k>\lambda n} a_{n,k}\underset{e}\to\sim \hat r_{+,\lambda}(1)^{-n},
$$
where
$
r_{+,\lambda}(1)$ is equal to $r(z_{0,\lambda})z_{0,\lambda}^\lambda $ if  $1<z_{0,\lambda}$
and to $r(1)$ otherwise.
Due to logarithmic convexity of $D$,
$$
z_{0,\lambda}>z_{0}\quad\text{and}\quad
r(z_{0,\lambda})z_{0,\lambda}^\lambda>r(z_{0})z_{0}^\lambda>r(z_0)z_0^\frac12,
$$
if $\lambda>\frac12.$
\qed\enddemo

\head  2. On the impact of Arnold inequalities.
\endhead

\subhead 2.1. Impact of the bound on the number of nonempty ovals
\endsubhead
Consider first the case of curves of degree $d$
with ${(d-1)(d-2)}/2+1$
connected components and denote by $L_d$ the number of the arrangements which
satisfy the Arnold bound
(5).
Encoding the arrangements by trees
we find that $L_d$ is the number of rooted unlabelled trees with
$n=N_d$
vertices and $\ge K_d$
leaves where $K_d=N_d-1-[\frac{d-1}2]([\frac{d-1}2]-1)$.
Recall that $N_d\sim \frac{d^2}2$ (see also (6)).

\proclaim{Proposition 8}
$$
    L_d\underset{e}\to\sim C_1^{\frac{d^2}2}.
$$
\endproclaim
\demo{Proof}
We apply Theorem 7. Since $C(\lambda)$ is continuous at
$\lambda=\frac12$, we find 
for any $\epsilon>0$ such $\delta>0$ that for any sufficiently big 
$n$ it holds
$$
(C_1+\epsilon)^{(1+\epsilon)n}\ge\sum_{k>(\frac12-\delta)n}a_{n,k}
\quad\text{and}\quad
\sum_{k>(\frac12+\delta)n}a_{n,k}\ge (C_1-\epsilon)^{(1-\epsilon)n}.
$$
It remains to put $n=N_d\sim\frac{d^2}2$ and to note that
for any sufficiently big $d$
$$
\sum_{k>(\frac12-\delta)n}a_{n,k}\ge
L_d\ge \sum_{k>(\frac12+\delta)n}a_{n,k}.\qed
$$
\enddemo

Now, consider the general case and denote, in accordance
with
the
Arnold bound on the number of empty ovals, by $L'_d$ the number of rooted unlabelled
trees with $n\le N_d$
vertices and $\ge n-[\frac{d-1}2]([\frac{d-1}2]-1)$
leaves.

\proclaim{Theorem 9}
$$
    L'_d\underset{e}\to\sim C_1^{\frac{d^2}2}.
$$
\endproclaim
\demo{Proof} In view of
(\eqOtter) and Proposition 8, it is sufficient to
prove that
$L'_d\le (\hat k^2-\hat k)T_{\hat k^2-\hat k}+k^2 L_d$
where $\hat k=[\frac{d-1}2]$
and $k=[\frac{d}2]$. Clearly, the first term bounds from above the
total number of trees with $n\le\hat k^2-\hat k$ vertices.
In the range $\hat k^2-\hat k<n\le N_d$ the number of the trees excluded
by the Arnold bound (\eqLbound) is increasing, from $0$ to $L_d$, when $n$ grows, since
$a_{n,m}\le a_{n+1,m+1}$  (to prove
such an inequality it is sufficient to add a leaf to a branch with a maximal number
of leaves). The coefficient $k^2$ before $L_d$ is due to
$$
N_d-1-[\frac{d-1}2]([\frac{d-1}2]-1)=k^2.\qed
$$\enddemo

\subhead 2.2. Auxiliary lemmas
\endsubhead

Let $v$ be a vertex of a tree $t$.
A {\it branch of $t$ at $v$} is a connected component of
the graph obtained from $t$ by removing $v$ and the
(open) edges adjacent to $v$.

\proclaim{ Lemma \lempnOne } Let $t$ be a tree with $N$ vertices.
Then there exists a vertex $v$ such that any branch of $t$ at $v$
has at most $N/2$ vertices.
\endproclaim

\demo{ Proof }
Suppose that any vertex has a branch with more than $N/2$ vertices.
Choose any vertex $v_1$ and
define the sequence of
vertices $v_1,v_2,\dots$ as follows.
Assume that $v_i$ is already defined.
Let $t_i$ be the branch of $t$ at $v_i$ which has more than $N/2$
vertices. Then $v_{i+1}$ is defined as the vertex of $t_i$ which
is nearest to $v_i$.
Moving from $v_1$ to $v_2$, then from $v_2$ to $v_3$ and so on, we
can never turn back. Indeed, if $v_{i+1}$ coincides with $v_{i-1}$ then
removing from $t$
the (open) edge connecting $v_i$ with $v_{i+1}$
we would obtain two
subtrees of $t$ each having more than $N/2$ vertices.
Since $t$ has no loops, this means that our sequence has no repeatings.
Contradiction.
\qed
\enddemo

\proclaim{ Lemma \lempnTwo }
Let $c_1\ge\dots\ge c_r\ge 0$
and
$|c|\le c_1+\dots+c_r$. Then there exist
$\eps_1,\dots,\eps_r\in\{\pm1\}$
such that $|(\eps_2 c_2+\dots+\eps_r c_r)-c|\le c_1$.
\endproclaim

\demo{ Proof }
Set $\eps_{k+1} = \sign_+(c-(\eps_2 c_2+\dots+\eps_k c_k))$ where
$\sign_+(x) = 1$ for $x\ge0$ and $\sign_+(x)= -1$ for $x<0$.
This means that we walk along the real axis starting
from the origin so that the absolute values of the steps are
successively $c_2,c_3,...$ and each step is directed
towards the point $c$.
Then $c'=
\eps_2 c_2+\dots+\eps_r c_r$
is the final point of our walk. It is easy to see that
$|c'-c| \le c_1$.
\qed
\enddemo

In accordance with the terminology coming from the geometry
of plane curves, let us say that a vertex of a rooted tree
$t$
is {\it even}
(resp. {\it odd})
if the minimal path relating it to the root consists of an
odd
(resp.
even) number of edges.
Let denote by $p(t)$ (resp.  $n(t)$)
the number of even (resp. odd) vertices, including the root, of $t$
and put $\chi(t) =
p(t)-n(t)$.

For example, the root is an
odd vertex, the vertices connected to the root by an edge are
even etc. Note, that when we change the root, $|\chi(t)|$
does not change.

We say that a rooted tree $t'$ is obtained from a rooted tree $t$
by {\it contracting an edge}
if $t'$ is obtained from $t$ by replacing
some edge with a single vertex $v$ (see Figure \figRemoveEdge).
If one of the ends of the edge
which we contracted was the root of $t$, then
$v$ is declared the root of $t'$.
This operation reduces the number of vertices and the edges by one.
The operation of {\it inserting an edge at $v$} is to be thought of
as an inverse operation.
When one of the ends of the inserted
edge is a leaf, this is called the {\it attachment of an edge}.

\midinsert
\epsfxsize 80mm
\centerline{\epsfbox{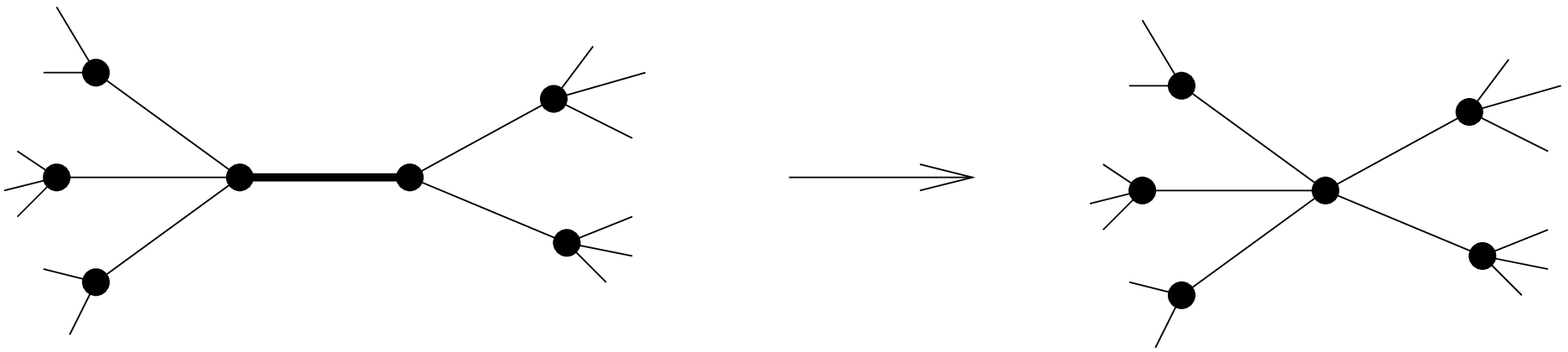}}
\botcaption{ Figure \figRemoveEdge. Edge removing }
\endcaption
\endinsert

\proclaim{ Lemma \lempnThree }
Let $t_0$ be a rooted tree with $N$ vertices and let $c$ be any integer
such that $|c| \le |\chi(t_0)|$.
Then there exists a sequence of rooted trees $t_1,\dots,t_k$
such that
\roster
\item $\chi(t_k)=c$;
\item $t_{i-1}$, $i=1,\dots,k$, is obtained from $t_i$
      by
      contracting an edge;
\item $k \le 3 + 3\log_2 N$.
\endroster
\endproclaim

\demo{ Proof }
Apply the induction by $N$.
The case $N=1$ is trivial.
Assume that the statement is true for any tree which has
less than $N>1$ vertices.
By Lemma \lempnOne, there exists a vertex $v$ such that
any branch of $t_0$ has at most $N/2$ vertices.
Let us denote the branches of $t$ at $v$ by $b_1,\dots,b_r$.
We
 choose the root of each branch at the vertex nearest to $v$.
Let $c_i=|\chi(b_i)|$ and $\delta_i=\sign\chi(b_i)$.
Let us number the branches so that
$c_1\ge c_2\ge\dots\ge c_r$.
By Lemma \lempnTwo, there exist $\eps_2,\dots,\eps_r\in\{\pm1\}$
such that $|c'-c|\le c_1$ where $c'=\eps_2 c_2+\dots+\eps_r c_r$.
By the induction hypothesis, we can insert $\le 3+3\log(N/2)=3\log N$
edges to $b_1$ so that $c_1^*=\chi(b_1^*)=c-c'$
for the resulting tree $b_1^*$. Let $t_0^*$ be the tree obtained from $t_0$
by replacing $b_1$ with $b_1^*$.

Let $t_1^*$ be obtained from $t_0^*$ by inserting an edge $e$ at $v$
so that $b_1^*$ and the branches $b_i, i\ge 1,$ with $\eps_i\delta_i
=\sign(c-c')$
are on one side of $e$
and the branches $b_i$ with $\eps_i\delta_i=
-\sign(c-c')$
are on
the other side. Then we have
$|\chi(t_1^*)|=|c_1^*+c'|=|c|$. Now,
we may return to counting $|\chi|$ with
respect to the initial root of $t_0$ and respective roots of $t_i, i\ge 1$.
If $\chi(t_1^*)=-c$, we attach an edge to the root, choose the
obtained leaf as the new root and then attach an edge to the new root.
\qed
\enddemo

\subhead 2.3. Impact of the bounds on the number of even and odd nonempty ovals
\endsubhead

Let us
denote by $A_d$ the number of rooted unlabelled
trees with $n\le N_d$ vertices
which satisfy the Arnold bounds (\eqAr).

\proclaim{Theorem 13}
$$
    A_d\underset{e}\to\sim C_1^{\frac{d^2}2}.
$$
\endproclaim
\demo{Proof}
If a tree with $n\le N_d$ vertices
satisfies the weak Arnold bound (\eqLbound), we apply to it, removing its leaves,
Lemma \lempnThree \, with $c=0$, and then put the leaves back,
getting thus
a tree with $n+3[\log_2 n]+3\le N_d+3[\log_2 N_d]+3\le N_{d+6}$ vertices which
satisfies the stronger Arnold bounds
(\eqAr). Therefore,
$$
    {L'_{d}\over A_{d+6}} \le
    \sum_{n=1}^{N_d}\binom {n+3[\log_2 n]+3} {3[\log_2 n]+3}
    \le N_d\binom {N_d+3[\log_2 N_d]+3} {3[\log_2 N_d]+3} = e^{o(N_d)}
$$
and the theorem follows now
from Theorem 9 and $A_d\le L'_d$.
\qed\enddemo


\head  Appendix. Limit distribution
\endhead

Let us consider $a_{n,k}/a_n(1)$ as a probability distribution
of a random variable $X_n$, i.e. $P(X_n=k)=a_{n,k}/a_n(1)$.
As is known, see f.e. \cite{\refDG}, the following central 
limit theorem holds:
this
random sequence $X_n$, once normalized, tends
to a normal distribution:
$$
P(a<\frac{X_n-mn}{\sigma\sqrt{n}}<b)\to\frac1{2\pi}\int_a^be^{-\frac{x^2}2}dx.
$$
where
$$
    m=-r'(1)/\alpha = 0.4381562356643746639684921638628797837055...
$$
and
$$
    \sigma^2 = {r'(1)^2\over\alpha^2} - {r'(1)+r''(1)\over\alpha} =
    0.150044811672846981980699640444640111071...
$$
In particular, this means that approximately $43.8\%$ of vertices of
a big random tree are leaves. It is worth mentioning that the fact
that the mean value of the number of leaves is $\sim mn, m=0.438156235664\dots$
was established
by Robinson and Schwenk \cite{\refRS} by the Polya-Otter method,
and its extension to the other moments was given by Schwenk in \cite{\refS}.

In view of the above limit theorem, it is natural to replace
$a_{n,k}$ by its approximation by the normal distribution
$$
    a_{n,k}^* = {a_n(1)\over\sigma\sqrt{2\pi}}
    \exp\left(-{(k-mn)^2\over2\sigma^2n}\right). 
$$
Then, we get
$$
    \sum_{k>n/2} a_{n,k}^* \underset{e}\to\sim
    \alpha^{-n} \exp\left(-{(1/2-m)^2 n\over 2\sigma^2}\right) = C_2^n
$$
where
$$
    C_2 = \alpha^{-1}\exp\left(-{(1/2-m)^2\over 2\sigma^2}\right) =
    2.91833301345955740149786987821329181193...
$$

We see that $C_2$ differs from $C_1$ in the fourth digit.
This is not a contradiction with the central limit theorem because
this
just means that the convergence to the normal distribution is
not good far from the center. It shows
that the central limit theorem is not sufficient for a search
of the rate of growth of $\sum_{k>n/2} a_{n,k}$.

To conclude, let us notice that
the
constants $r'(1)$ and $r''(1)$ (needed to find
$m$ and $\sigma^2$) can be computed much faster than the constants
$z_0$ and $x_0$ from Section 2
because the double
summation over $n,k$ may be replaced with the single summation by
use of the following recurrent formulas for the coefficients of
of the series $T_z(x,1)$ and $T_{zz}(x,1)$.
Similarly to (\eqRR), one can obtain
$$
    a'_{n+1}(1)=\sum_{j=1}^n a'_j(1) \sum_{k=1}^{[n/j]} a_{n+1-kj}(1)
$$
$$
\split
    a''_{n+1}(1)=&\sum_{j=1}^n\Big\{
            a'_j(1)\Big(\sum_{k=1}^{[(n-1)/j]} a'_{n+1-kj}(1)\Big)
\\
            &+ a'_j(1)\Big(\sum_{k=1}^{[n/j]}(k-1)a_{n+1-kj}(1)\Big)
            + a''_j(1)\Big(\sum_{k=1}^{[n/k]}ka_{n+1-kj}(1)\Big)\Big\}
\endsplit
$$


\enddocument